\documentclass[11pt]{amsart}

\usepackage{verbatim, amscd, epsfig,amsmath,amsthm,amstext,amssymb,a4,hyperref}
\theoremstyle{plain}
\newtheorem{theorem}{Theorem}[section]
\newtheorem{lemma}[theorem]{Lemma}
\newtheorem{cor}[theorem]{Corollary}
\newtheorem{prop}[theorem]{Proposition}

\theoremstyle{definition}
\newtheorem{definition}[theorem]{Definition}
\newtheorem{example}[theorem]{Example}
\newtheorem{examples}[theorem]{Examples}

\theoremstyle{remark}
\newtheorem{rem}[theorem]{Remark}

\numberwithin{equation}{section}

\DeclareMathOperator{\Hom}{Hom}
\DeclareMathOperator{\End}{End}

\DeclareMathOperator{\tr}{tr}
\DeclareMathOperator{\id}{id}
\DeclareMathOperator{\dbar}{{\bar{\partial}}}
\DeclareMathOperator{\delbar}{{\bar{\partial}}}
\DeclareMathOperator{\rank}{rank}
\DeclareMathOperator{\ord}{ord}

\newcommand{\trivial}[1]{\underline{\H}^{#1}}

\newcommand{\thedef}[1]{\emph{#1}}

\newcommand{\invers}{^{-1}}
\newcommand{\einhalb}{\frac 12}
\newcommand{\einviertel}{\frac 14}

\newcommand{\Zz}{{\mathcal Z}}

\renewcommand{\Im}{{\rm Im}\,}
\newcommand{\del}{\partial}
\newcommand{\R}{\mathbb{ R}}
\newcommand{\C}{\mathbb{ C}}

\renewcommand{\H}{\mathbb{ H}}

\newcommand{\N}{\mathbb{ N}}
\renewcommand{\P}{\mathbb{ P}}
\newcommand{\HP}{\H\P}
\newcommand{\CP}{\C\P}

\begin{document}

\title
{Willmore spheres in quaternionic projective space}
\author{K. Leschke}
\address{
Katrin Leschke\\
Institut f\"ur Mathematik\\
Lehrstuhl f\"ur Analysis und Geometrie\\
Universit\"at Augsburg\\
D-86135 Augsburg
}
\email{katrin.leschke@math.uni-augsburg.de}
\thanks{MSC--class:  53Axx, 53Cxx, 30Fxx \\
Partially supported by SFB 288  and   by  NSF-grant DMS-9626804}
\begin{abstract}

The Willmore energy for Frenet curves in quaternionic projective space
$\HP^n$ is the generalization of the Willmore functional for
immersions into $S^4$.  Critical points of the Willmore energy are
called Willmore curves in $\HP^n$.

Using a B\"acklund transformation on Willmore curves, we generalize
Bryant's result on Willmore spheres in 3--space: a Willmore sphere in
$\HP^n$ has integer Willmore energy, and is given by complex
holomorphic data.
\end{abstract}


\maketitle


\section{Introduction} 
This work is part of a project where ``quaternionified'' complex
analysis is used to study old and new questions in surface theory.
The first accounts of this program are presented in \cite{icm},
\cite{coimbra} and \cite{Klassiker}. An important feature in the
quaternionic setup is that conformal maps from a Riemann surface into
$S^4 = \HP^1$ play the role of the meromorphic functions in complex
analysis.  More generally, if we consider holomorphic curves $f: M \to
\HP^n$ then the components of $f$ are branched conformal immersions
into $S^4 = \HP^1$. Thus we can think of a holomorphic curve in
$\HP^n$ as a family of branched conformal immersions into $S^4$.

Basic constructions of complex Riemann surface theory, such as
holomorphic line bundles, the Kodaira embedding, the Pl\"ucker
relations and the Riemann--Roch theorem, carry over to the
quaternionic setting. There is an important new invariant of the
quaternionic holomorphic theory distinguishing it from its complex
counterpart: the Willmore energy is defined for holomorphic curves in
$\HP^n$.  For immersions $f: M \to S^4$ into $S^4= \HP^1$ we obtain
the classical Willmore functional $\int_M(|H|^2 - K - K^\perp)|df|^2$
where $H$ is the mean curvature vector, $K$ the Gaussian curvature,
and $K^\perp$ the curvature of the normal bundle. Willmore surfaces
$f: M \to S^4$ are critical points of the Willmore functional
\cite{willmore_book}. It is a well--known fact, that Willmore surfaces
in $S^4$ are characterized by the harmonicity of the conformal Gauss
map.

To generalize the notion of Willmore surfaces to the case of a
holomorphic curve $f: M \to \HP^n$, we use the the analogue of the
conformal Gauss map of an immersion into $S^4$, the so--called
\emph{canonical complex structure} of $f$.  In general, the canonical
complex structure exists only away from a discrete set of $M$. In view
of the relation between the Willmore condition and harmonicity, we
restrict to the case of Frenet curves.  These are holomorphic curves
for which the canonical complex structure exists smoothly on $M$.  A
Frenet curve in $\HP^n$ is called \emph{Willmore} if it is a critical
point of the Willmore energy under compactly supported variations by
Frenet curves.

Similar to the $\delbar$ and $\partial$ transforms of harmonic maps
into complex projective space \cite{wolfson}, we define a B\"acklund
transform of a Willmore curve by using the $(1,0)$--part of the
derivative of the canonical complex structure. This generalizes the
B\"acklund transformation in \cite{coimbra} for Willmore surfaces in
$\HP^1$ to Willmore curves in $\HP^n$.

The B\"acklund transform $\tilde f$ of a Willmore sphere $f: S^2\to
\HP^n$ is a Willmore sphere in $\HP^k$, $k<n$, or a constant point in
$\HP^n$.  More precisely, the B\"acklund transform $\tilde f: S^2
\to\HP^k$ is given by a twistor projection of a holomorphic curve $g:
S^2 \to \CP^{2k+1}$.  Using the special form of the B\"acklund
transform, we prove a generalization of the results of Bryant
\cite{bryant} and Ejiri, \cite{ejiri}: every Willmore sphere
$f:S^2\to\HP^n$ is either a minimal surface in $\R^4$ with planar ends
or is given by a twistor projection of a holomorphic curve in some
complex projective space.  Moreover, a Willmore sphere has Willmore
energy in $4\pi \N$.



\section{Holomorphic curves and holomorphic bundles}

We set up some basic notation used throughout the paper. For more
details of the underlying quaternionic theory, we refer to
\cite{coimbra}, \cite{Klassiker} and \cite{icm}.

We view a Riemann surface $M$ as a $2$--dimensional, real manifold
with an endomorphism field $J\in\Gamma(\End(TM))$ satisfying $J^2=-1$.
If $V$ is a vector bundle over $M$, we denote the space of $V$ valued
quaternionic $k$--forms by $\Omega^{k}(V)$.
If $\omega\in\Omega^{1}(V)$, we set
\begin{equation*}
*\omega:=\omega\circ J\,.
\end{equation*}
Moreover, we will identify $\omega\in\Omega^2(V)$ with the induced
quadratic form $\omega(X):=\omega(X,JX)$. In particular for pairings
$V_1\times V_2\to V_3$ we will identify
\[ \omega\wedge \eta = \omega*\eta - *\omega\eta, \ 
\omega\in\Omega^1(V_1), \eta\in\Omega^1(V_2),\] where the wedge
product is defined over the pairing.

Most of the vector bundles occurring will be {\em quaternionic} vector
bundles, i.e., the fibers are quaternionic vector spaces and the local
trivializations are quaternionic linear on each fiber.  We adopt the
convention that all quaternionic vector spaces are {\em right} vector
spaces.  A {\em quaternionic connection} on a quaternionic vector
bundle satisfies the usual Leibniz rule over quaternionic valued
functions.

If $V_1$ and $V_2$ are quaternionic vector bundles, we denote by
$\text{Hom}(V_1,V_2)$ the bundle of quaternionic linear homomorphisms.
As usual, $\text{End}(V)=\text{Hom}(V,V)$ denotes the quaternionic
linear endomorphisms.  Notice that $\text{Hom}(V_1,V_2)$ is {\em not}
a quaternionic bundle.

Let $V$ be a quaternionic vector bundle with complex structure
$S\in\Gamma(\End(V))$, $S^2 = -1$. We phrase this as $(V,S)$ is a
complex quaternionic vector bundle.

Given a connection $\nabla$ on $V$, we can decompose $\nabla = \nabla'
+ \nabla''$ into $(1,0)$ and $(0,1)$ parts with respect to $S$, where
 \[
\nabla' := \einhalb(\nabla - S *\nabla) \ 
 \text{ and } \nabla'' := \einhalb(\nabla + S *\nabla)\,.
\]
 
 Let $\Hom_\pm(V,W) =\{ B\in\Hom(V,W) \mid S_WB = \pm BS_V\}$ where
 $(V,S_V)$ and $(W,S_W)$ are complex quaternionic vector bundles. We
 denote by $B_\pm = \einhalb(B \mp S_WBS_V)\in\Hom_\pm(V,W)$ the
 $\pm$--part of $B\in\Hom(V,W)$. We can decompose $\nabla'$ and
 $\nabla''$ further into
\[
\nabla'' = \dbar  + Q, \ \nabla' = \partial + A\,
\]
where
 \[
\dbar S = S\dbar, \ \partial S = S\partial
\]
 and 
\[
Q = \nabla''_-\in\Gamma(\bar K \End_-(V)), \ A=
\nabla'_-\in\Gamma(K\End_-(V))\,.
\]
Here we denote by
\[
 KE :=
\{\omega\in \Lambda^1(TM) \otimes E\mid *\omega = S\omega\}
\]
 and 
\[
\bar KE :=
\{\omega\in \Lambda^1(TM) \otimes E\mid *\omega = -S\omega\}
\]
the tensor product of the canonical and anti-canonical bundle with a
complex vector bundle $(E,S)$.  Notice that
\begin{equation}
\label{eq:nablaS} \nabla S = [\nabla, S] = [Q+A,S] = 2(*Q -*A)\,.
\end{equation}
 In particular, if $\nabla$ is flat, then 
\begin{equation}
\label{eq:dnabla*A}
 d^\nabla *A = d^\nabla *Q 
\end{equation}
and
\begin{equation}
\label{eq:*A}
 4*A = S*\nabla S - \nabla S, \quad 4*Q = S*\nabla S + \nabla S\,.
\end{equation}

\quad

We denote by $V$ the trivial $\H^{n+1}$--bundle over
$M$.  Let $\Sigma\to G_{k+1}(V)$ be the tautological $(k+1)$--plane
bundle whose fiber over $V_k\in G_{k+1}(\H^{n+1})$ is
$\Sigma_{V_k}=V_k\subset V$.  A map $f:M\to G_{k+1}(V)$ can be
identified with a rank $k+1$ subbundle $V_k\subset V$ via
$V_k=f^{*}\Sigma$, i.e., $(V_k)_p=\Sigma_{f(p)}=f(p)$ for $p\in M$.
From now on, we will make no distinction between a map $f$ into the
Grassmannian $G_{k+1}(V)$ and the corresponding subbundle $V_k\subset
V$.

The derivative of $V_k\subset V$ is given by the
$\text{Hom}(V_k,V/V_k)$ valued $1$--form
\begin{equation}\label{eq:derivative}
\delta=\pi_{V_k} \nabla|_{V_k}\,,
\end{equation}
where $\pi_{V_k}:V\to V/V_k$ is the canonical projection. Under the
identification
$TG_{k+1}(V)=\text{Hom}(\Sigma,V/\Sigma)$ the
$1$--form $\delta$ is the derivative $df$ of $f:M\to G_{k+1}(V)$.\\

\begin{definition}
\label{def:Frenetflag}
Let $V=\trivial{n+1}$ be the trivial quaternionic $n+1$--plane bundle
over a Riemann surface $M$.
\begin{enumerate}
\item A rank $k+1$ subbundle $V_k\subset V$ is a {\em holomorphic
    curve} in $V$ if there exists a complex structure
  $J\in\Gamma(\text{End}(V_k))$, $J^2=-1$, such that
\begin{equation*}
*\delta=\delta\,J\,.
\end{equation*}
\item The \emph{Frenet flag} of a holomorphic curve $f: M \to\HP^n$ is
  a full flag 
\[
 L =f^*\Sigma = V_0 \subset V_1 \subset \ldots \subset V_{n-1}
  \subset V_n = V
\]
of quaternionic subbundles of rank $V_k = k+1$ together with complex
structures $J_k$ on the quotient bundles $V_k/V_{k-1}$, $J_0 = J$,
such that
\begin{enumerate}
\item $\nabla\Gamma(V_k)\subset \Omega^1(V_{k+1})$.
\item The derivatives $\delta_k = \pi_{V_k}\nabla: V_k/V_{k-1} \to
  T^*M\otimes V_{k+1}/V_k$ satisfy
\[ 
*\delta_k = J_{k+1}\delta_k = \delta_k J_k\,.
\]
\end{enumerate}
\end{enumerate}
\end{definition}

\quad

\begin{example} Consider the line bundle $L= f^*\Sigma \subset
  \underline{\H}^2$ induced by $f: M \to S^4 = \HP^1$. The
  corresponding flag is $L \subset V_1 = \trivial{2}$. The line bundle
  $L$ is a holomorphic curve if there exists a complex structure $J$
  on $L$ such that $*\delta = \delta J$. We obtain an equivalent
  condition written in terms of $f$: $L$ is a holomorphic curve if and
  only if there exists $J\in \Gamma(\End(L))$ such that $*df = df
  \circ J$, i.e. $f: M \to \HP^1$ is a branched conformal immersion
  (for a detailed development of conformal surface theory using
  quaternionic valued functions see \cite{coimbra}).  If the
  derivative of $L$ has no zeros, then $f$ is an immersion. In this
  case, we can define the conformal Gauss map which can be identified
  with a complex structure $S$ on $\trivial{2}$ which induces $J$ on
  $L$.  Moreover, $S$ satisfies $*\delta = S\delta = \delta S$ and the
  second order tangency condition $Q|_L=0$.  The immersion $f$ is
  Willmore if and only if the conformal Gauss map is harmonic
  \cite{ejiri}, \cite{rigoli}, which is equivalent 
  \cite[Thm.\ 3]{coimbra} to 
  \[
  d^\nabla * A = 0\,.
\]
\end{example}

Recall that a complex structure $S\in\Gamma(\End(V))$ is called
\emph{adapted} to the Frenet flag $L=V_0 \subset V_1 \subset \ldots
\subset V$ of a holomorphic curve $f: M \to\HP^n$  if $S$ induces the
complex structures given by the Frenet flag, i.e., if 
\[
 *\delta_k = S\delta_k = \delta_k S\,
\]
for $k=0, \ldots, n-1$.  The analogue of the conformal Gauss map of a
conformal immersion $f: M \to S^4$ is  an adapted complex structure
which satisfies a certain second order condition.

\begin{definition} 
  Let $f: M \to\HP^n$ be a holomorphic curve and $L=V_0 \subset V_1
  \subset \ldots \subset V$ be the Frenet flag of $L$.  The unique
  \emph{adapted} complex structure $S\in\Gamma(\End(V))$ with
\[
 Q|_{V_{n-1}} = 0 \ \text{or, equivalently } A(V) \subset TM^* \otimes
 L\,,
\]
is called the \emph{canonical complex structure} of $L$.  Here $\nabla
= \dbar + Q + \partial + A$ is the decomposition of $\nabla$ with
respect to $S$.
\end{definition}
In general, the Frenet flag and the canonical complex structure of a
holomorphic curve $f: M \to\HP^n$ only exist away
from a discrete set $D$, the Weierstrass points of $L$.  These are the
zeros of the derivatives $\delta_k$ of the flag bundles $V_k$. In case
of a holomorphic curve $f$ in $\HP^1$ the Weierstrass points are the
branch points of the map $f: M \to S^4$.

Whereas the Frenet flag of a holomorphic curve always extends
continuously into the Weierstrass points \cite[Lemma 4.10]{Klassiker},
the canonical complex structure may become singular as the following
example \cite{paule} shows: If $f: M \to\HP^1$ is the twistor
projection of a complex holomorphic curve $h: M \to\CP^3$ then the
canonical complex structure is given by the tangent line $W_1\subset
V$ of $h$, namely $S|_{W_1} = i$ and $S|_{W_1j} = -i$. But the tangent
$W_1\subset V$ of $h$ can become quaternionic, i.e., $W_1 = W_1 j$ at
some $p\in M$. In this case the canonical complex structure $S$
degenerates to a point at $p\in M$ and thus $S$ cannot be extended
into $p\in M$.  To avoid these difficulties, we will only consider
holomorphic curves $f: M \to\HP^n$ which have a smooth canonical
complex structure. For conformal maps $f: M \to \HP^1$ this means that
the mean curvature sphere congruence extends smoothly across the
branch points.

\begin{definition}
  A \emph{Frenet curve} $f: M \to\HP^n$ is a holomorphic curve which
  has a smooth canonical complex structure on $M$.
\end{definition}

\begin{rem}
  Note that this definition is more general then the one given in
  \cite{Klassiker}. In contrast to \cite[Def.\ 4.3]{Klassiker} we
  allow $f$ to have Weierstrass points. However, the smoothness of the
  canonical complex structure guarantees, as we will see below, the
  existence of the Frenet flag on $M$.
\end{rem}

\begin{example} A trivial example for a Frenet curve is an \emph{unramified}
  curve $f: M \to\HP^n$, i.e., a holomorphic curve without Weierstrass
  points. In the case $n=1$ an unramified curve is a conformal
  immersion $f: M \to S^4$.
  
  In what follows, the example of the dual curve of a Frenet curve
  also will play an important role. For any subbundle $V_k$ of $V$ let
 \[
V_k^\perp := \{\alpha\in V\invers \mid \, <\alpha,
\psi> = 0  \text{ for all } \psi\in V_k\}\,,
\]
where $V\invers$ is the dual bundle of $V$.  The \thedef{dual curve}
$L^*$ of a Frenet curve $L\subset V$ with Frenet flag $L=V_0\subset
V_1\subset\ldots\subset V$ is the holomorphic curve in $V\invers$
defined by
\[
L^* := V_{n-1}^\perp\,.
\] 
The Frenet flag of the dual curve $f^*$ is given by
\begin{equation}
\label{eq:dualflag}
V^*_k =V^\perp_{n-1-k}\,.
\end{equation}
An adapted complex structure $S$ of $L$ induces an adapted
complex structure on $L^*$ by the dual map $S^*$.  If $\nabla =
\dbar + \partial + A + Q$ is a connection on $V$ then the induced dual
connection $\nabla^*$ on $V\invers$ is decomposed as
\[
\nabla^*= \hat \nabla^* + Q^\dagger + A^\dagger\,,
\]
where 
\begin{equation}
\label{eq:Adagger}
 A^\dagger = -Q^* \in\Gamma(K\End_-(V\invers))  
\quad \text{ and } \quad 
Q^\dagger = -A^* \in\Gamma(\bar K\End_-(V\invers))\,.
\end{equation}
The decomposition of the dual connection $\hat \nabla^* = \dbar^* +
\partial^*$ is given by
\begin{equation}
\label{eq:dbardual}
 <\dbar^* \alpha, \psi > + <\alpha,\dbar\psi> = \einhalb(d<\alpha,\psi> +
*d<\alpha, S\psi>)
\end{equation}
(and a corresponding equation for $\partial^*$).  If $S$ is the
canonical complex structure of $L$ then
\[
 Q^\dagger|_{V_{n-1}^*} = - A^*|_{L^\perp} = 0\,,
\]
and $S^*$ is the canonical complex structure of the dual curve $L^*$.
In particular, the dual curve of a Frenet curve is Frenet.  In the
case of a Frenet curve $f: M \to S^4$, the dual curve is given by
the antipodal map since $L^* = V_{n-1}^\perp = L^\perp$.
\end{example}

A holomorphic curve $f: M \to\HP^n$ induces a holomorphic structure
$D$ on the dual bundle $L\invers$ of $L=f^*\Sigma$. Recall the
definition of a (quaternionic) holomorphic structure:

\begin{definition}[see \cite{icm}]
  A \thedef{holomorphic structure} on a complex quaternionic vector
  bundle $(V,S)$ is a real linear map $D: \Gamma(V) \to \Gamma(\bar K
  V)$ satisfying 
\[
 D(\psi \lambda) = (D\psi)\lambda + \einhalb(\psi d\lambda + S\psi
*d\lambda), \ \lambda: M \to \H\,.
\] 
We denote by $H^0(V) = \ker D \subset \Gamma(V)$ the space of
holomorphic sections and call $(V,S,D)$ a \thedef{holomorphic
  quaternionic vector bundle}.  A subbundle $W\subset V$ is called a
\thedef{holomorphic subbundle} of $(V,S,D)$ if $\Gamma(W)$ is $D$
stable. In this case $(W,S|_W,D|_W)$ is a holomorphic quaternionic
vector bundle.
\end{definition}

We decompose a holomorphic structure into 
\[
D = \dbar + Q\,,
\]
where $\dbar S = S\dbar$ and $Q\in\Gamma(\bar K \End_-(V))$. The case
when $Q=0$ gives the usual theory of (doubles) of complex holomorphic
vector bundles. We call $D = \dbar$ a \emph{complex holomorphic
  structure}.

\begin{examples}
\label{ex:holo_structure}
\begin{enumerate}
\item 
If $\nabla $ is a connection on $V$ then $\nabla''$ is a
holomorphic structure on $(V,S)$.
\item We denote by $\bar V$ the complex vector bundle $(V,-S)$ of the
  complex quaternionic vector bundle $(V,S)$. Then $\nabla'$ is the
  $(0,1)$--part of the connection $\nabla$ with respect to $-S$ and
  $(\bar V, \nabla')$ is a holomorphic vector bundle. We call
  $\nabla'$ an \emph{antiholomorphic structure} on $(V,S)$.
\item 
If $V_k$ is an $S$ stable subbundle of $(V,S)$ which is also a
holomorphic curve with respect to $S|_{V_k}$, i.e.,
\[ 
*\delta = \delta S, \ \text{ where } \delta= \pi_{V_k} \nabla|_{V_k} :
 V_k \to V/V_k\,,
\]
then $V_k$ is in general \emph{not} a holomorphic subbundle of $V$
with respect to $\nabla''$. The condition for $V_k$ being a
holomorphic subbundle is exactly
\begin{equation*}
*\delta = S\delta\,,
\end{equation*}
which is equivalent to the condition that $V_k^\perp$ is a holomorphic
curve.

In particular, if $f: M \to\HP^N$ is a Frenet curve and $S$ is an
adapted complex structure of $f$ on $V$, then $S$ stabilizes the flag
spaces $V_k$ and satisfies $*\delta_k = S\delta_k = \delta_k S$. Thus
the $V_k$'s are holomorphic subbundles of $V$ with respect to
$\nabla''$ and holomorphic curves with respect to $S$.

\item However, the holomorphic curve $f: M \to\HP^n$ induces a
  canonical holomorphic structure on the dual bundle $L\invers$ of
  $L$: it is given by the requirement that restrictions of forms
  $\alpha\in(\H^{n+1})^*$ to $L$ are holomorphic sections. The space
  $\{\alpha|_L\mid \alpha\in (\H^{n+1})^*\}\subset H^0(L\invers)$ is a
  basepoint free linear system \cite[Sec.\ 2.6]{Klassiker}.
  
  In fact, the correspondence between holomorphic curves and basepoint
  free linear systems $H\subset H^0(L\invers)$ of a holomorphic line
  bundle $L\invers$ is one to one: by the \emph{Kodaira
    correspondence} one can consider the bundle $L$ as a subbundle of
  $V = H\invers$, see \cite[Thm.\ 2.8]{Klassiker}.
\end{enumerate}
\end{examples}

\begin{lemma}[see {\cite[Lemma 2.4]{coimbra}}]
  Let $S$ be a complex structure on $V$ and $\nabla = \dbar +
  \partial + Q + A$ the decomposition of the trivial connection on $V$
  with respect to $S$.  Then
\begin{equation}\label{eq:r2} R^{\dbar + \partial} = -(Q\wedge Q + A\wedge A)
  =  2S(A^2 - Q^2)
\end{equation}
and, for $Z\in H^0(TM)$,
\begin{equation}
  \label{eq:r1}
 R^{\dbar + \partial}_{Z,JZ} = 2S(\dbar_Z\partial_Z - \partial_Z\dbar_Z)\,.
\end{equation}
\end{lemma}

If $V_1$ and $V_2$ are two complex holomorphic vector bundles with
complex holomorphic structures $\dbar_k$, then $\Hom_+(V_1,V_2)$
inherits a complex holomorphic structure $\dbar $ via
\begin{equation}
\label{eq:dbar}
 (\dbar A)\psi := \dbar_2 (A\psi) - A(\dbar_1 \psi)\,.
\end{equation}
The usual tensor product construction for complex holomorphic
structures induces a complex holomorphic structure on
$K\Hom_+(V_1,V_2)$.

\begin{lemma}
\label{lem:holbundel} 
Let $S$ be a complex structure on $V$ and let $V_k\subset V$ be an $S$
stable subbundle.  Then the following statements are equivalent:
\begin{enumerate}
\item $V_k$
  is $A, Q$ and $\dbar$ stable.
\item $*\delta_k = S \delta_k = \delta_k S$ where $ \delta_k = \pi_{V_k}
\nabla|_{V_k}: V_k \to V/{V_k}.$
\end{enumerate}
In this case
\[ 
\delta_k = \pi_{V_k}  \partial \in H^0(K\Hom_+(V_k,V/V_k))
\]
is a holomorphic section. Here the  holomorphic structure on $V_k$ is
$\delbar$ and the holomorphic structure on $V/V_k$ is defined by
\[
\dbar \pi_{V_k} = \pi_{V_k} \dbar\,.
\]
\end{lemma}
\begin{proof} By a type consideration, a vector bundle is $\dbar+A +Q$
  stable if and only if it is stable under $\dbar$, $A$ and $Q$. Since
  \[
\delta_k = \pi_{V_k}(\del + \dbar  + Q + A)\,,
\]
we see that $V_k$ is $\dbar + Q + A$ stable if and only if $\delta_k
=\pi_{V_k}\del$. But this is equivalent to $*\delta_k = S\delta_k =
\delta_kS$ again by type considerations $*A = -AS, *Q = - SQ$, and
$*\delbar = -S\delbar$.
  
  Since $\dbar$ maps sections of $V_k$ to one-forms in $V_k$, we can
  define a holomorphic structure $\dbar$ on $V/V_k$ by  
\[
\dbar \pi_{V_k}= \pi_{V_k} \dbar\,.
\]
By (\ref{eq:r2}) we see that $R^{\hat\nabla}$ stabilizes $V_k$ and we
obtain for any local holomorphic sections $\psi\in H^0(V_k)$ and $Z\in
H^0(TM)$:
\begin{eqnarray*}
\dbar_Z(\delta_k(Z,\psi)) &=& \dbar_Z(\pi_{V_k}\partial_Z\psi)
 = \pi_{V_k}(\dbar_Z(\partial_Z\psi)) 
\stackrel{\eqref{eq:r1}}=  \pi_{V_k}(\partial_Z\dbar_Z \psi) = 
0\,.
\end{eqnarray*}
Thus, $\delta_k$ is holomorphic because it maps holomorphic sections
$\psi\in H^0(V_k)$ and $Z\in H^0(TM)$ to a holomorphic section
$\delta_k(Z,\psi)\in H^0(V/V_k)$.
\end{proof}

Since the flag derivatives of a Frenet curve are holomorphic sections
by the previous lemma, we see that a Frenet curve has a smooth Frenet
flag.

\begin{cor}[see 
  \cite{osculates}] The Frenet flag of a Frenet curve $f: M \to\HP^n$
  is smooth on $M$.
\end{cor}
\begin{proof}
  The canonical complex structure $S$ of the Frenet curve $f$ exists
  smoothly on $M$. In particular, $*\delta_0 = S \delta_0 = \delta_0
  S$ and the previous Lemma shows that $\delta_0\in
  H^0(K\Hom_+(L,V/L))$ is a holomorphic section. Therefore, the image
  of $\delta_0$ defines a smooth subbundle $V_1\subset V$.  Proceeding
  inductively, the Frenet flag exists smoothly on $M$.
\end{proof}

Using again a type argument as in the proof of Lemma
\ref{lem:holbundel}, we derive a criterion to decide whether a given
complex structure is the canonical complex structure of a Frenet
curve:
\begin{lemma}
\label{lem:QA=0}
Let $S$ be a complex structure on $V$. Assume that $L\subset V$ is a
holomorphic curve with respect to $S$ and that $L$ has a Frenet flag
$L\subset V_1 \subset \ldots \subset V_n = V$. Denote the derivative
of $V_k$ by $\delta_k$ and define 
\[
\delta^k:= \delta_{k-1} \circ
\ldots \circ \delta_0
\]
for $ k=1,\ldots, n$, where $\delta^0 = \id|_L$ denotes the identity
map of $V$ restricted to $L$. With the usual decomposition $\nabla =
\dbar +\partial + Q +A$ of $\nabla$ with respect to the complex
structure $S$, the following are equivalent
\begin{enumerate}
\item
$Q\delta^k = 0$ for all $k=0,\ldots, n-1$.
\item $L$ is a Frenet curve and $S$ is the canonical complex structure
  of $L$.
\end{enumerate}
\end{lemma}
\begin{proof}
  Since the canonical complex structure $S$ of a Frenet curve
  satisfies $Q|_{V_{n-1}} = 0$ and $\Im \delta_j \subset V_{j+1}/V_{j}
  \subset V_{n-1}/V_{j}$, we get $Q|_{\Im \delta_j}= 0$ for all
  $j=0,\ldots,n-2$.

  For the converse, observe that $Q\delta^i = 0$ for all $i\le k$
  implies that $V_{k}$ is $Q$ stable, i.e., the derivative of $V_k$ is
  given by 
\begin{equation}
\label{eq:deltak}
\delta_k = \pi_{V_k}(\partial + \dbar + A)|_{V_k}\,.
\end{equation}
We proceed by induction.  Since $L$ is a holomorphic curve with
respect to $S$, we have $*\delta_0 = \delta_0 S$. But then
(\ref{eq:deltak}) shows that $V_0=L$ is $A$ and $\dbar$ stable and we
get 
\[
*\delta_0=S\delta_0
\]
since $\delta_0 = \pi_L(\partial)$. In particular, the complex
structure of the Frenet flag, see Definition \ref{def:Frenetflag}, is
given by $J_1 = S$ on $V_1/L$, and $*\delta_1 = \delta_1S$.
Proceeding inductively, we see $J_k = S$ on ${V_k/V_{k-1}}$ and
$\delta_k = \pi_{V_k}\partial|_{V_k}$.  Hence $S$ is adapted to the
flag, and satisfies $Q|_{V_{n-1}}= 0$.
\end{proof}

We finish this section by a fact on the Hopf fields $A$ and $Q$ of a
Frenet curve which will allow later on to describe Willmore curves
with vanishing Willmore energy.

\begin{lemma}\label{l:AaufLnull}
  Let $f: M \to\HP^n$ be a holomorphic curve and $U\subset M$ an open
  subset of $M$ so that both the canonical complex structure $S$ and
  the Frenet flag $L\subset V_1 \subset \ldots \subset V_{n-1}\subset
  V$ are smooth on $U$. Then
\begin{enumerate}
\item If the restriction  of $A$ to $L$ vanishes on $U$ then
  $A$ vanishes on $U$.
\item If $Q$ has image in $V_{n-1}$ on $U$ then $Q$ vanishes on $U$.
\end{enumerate}
\end{lemma}

\begin{proof}
  Note first, that $Q|_{V_{n-1}}=0$ implies for
  $\varphi\in\Gamma(V_{n-1})$ that $(d^\nabla*A)\varphi=
  (d^\nabla*Q)\varphi = d^\nabla(*Q\varphi) + *Q\wedge\nabla\psi = *Q
  \wedge\delta_{n-1}\varphi = 0$ by type.
  
  If we assume that $A|_{V_k} = 0$ on $U$ for some
  $k\in\{0,\ldots,n-1\}$, then for  $X\in\Gamma(TU)$ and for $
  \psi\in\Gamma(V_{k}|_U)$ we obtain
\begin{eqnarray*}
 2A_X(\delta_k)_X\psi& =& (*A\wedge \delta_k)_{X,JX}\psi 
= (d^\nabla*A)_{X,JX}\psi - d^\nabla(*A\psi)_{X,JX} = 0\,,
\end{eqnarray*}
and hence $A|_{V_{k+1}} = 0$ on $U$.  By induction we see that
$A|_L=0$ on $U$ implies $A= 0$ on $U$.

Moreover, $\Im Q \subset \ker Q$ if and only if $\Im A^\dagger \subset
\ker A^\dagger$ because $Q^*=-A^\dagger$, 
\[
 \ker A^\dagger =
\ker Q^* = (\Im Q)^\perp\,,
\quad \text{
and} \quad
\Im A^\dagger = \Im Q^* = (\ker Q)^\perp\,.
\]
Since $\Im A^\dagger = L^*$, we see $A^\dagger|_{L^*} =0$ and conclude
$A^\dagger=0$ on $U$ by part \thetag{1}. Therefore, $Q=
-(A^\dagger)^*=0$ on $U$.
\end{proof}



\section{Willmore curves in $\HP^n$}

The Willmore functional of an immersion to $S^4$ can be generalized to
the Willmore energy of a holomorphic curve \cite{Klassiker}. We define
Willmore curves in $\HP^n$ as Frenet curves which are critical points
of the Willmore functional under compactly supported variations by
Frenet curves, \cite{osculates}. As in the case of Willmore surfaces
in $\R^3$ the Willmore condition is related to harmonicity: the
canonical complex structure of a Willmore curve in $\HP^n$ is
harmonic.

Recall that a holomorphic curve $f: M \to\HP^n$ induces a holomorphic
structure on the dual $L\invers$ of the line bundle $L=f^*\Sigma$.

\begin{definition} Let  $f: M \to\HP^n$ be a holomorphic curve from  a
  compact Riemann surface $M$ into quaternionic projective space.  The
  \thedef{Willmore energy} of $f$ is given by
\begin{equation}
\label{eq:Wenergy}
 W(f) := 2 \int_M < Q_{L\invers}\wedge *Q_{L\invers}>, 
\end{equation}
where $Q_{L\invers}$ is given by the holomorphic structure $D = \dbar
+ Q_{L\invers}$ on $L\invers$, and $<B> := \einviertel \tr_{\R} B$ for
$B\in\End(V)$.
\end{definition}
\begin{rem}
  Note that the definition of the Willmore energy is invariant under
  projective transformations of $f$. In the case $n=1$, we obtain
  \cite[Prop.\ 13]{coimbra} the usual Willmore functional
\[
  W(f) = \int_M (|H|^2 - K - K^\perp)|df|^2\,,
\]
of an immersion $f: M \to S^4$. Here $H$ is the mean curvature vector
of $f$, $K$ the Gaussian curvature, and $K^\perp$ the curvature of the
normal bundle
\end{rem}
From now on, $M$ will always denote a compact Riemann surface.

For a Frenet curve the Willmore energy can be computed in terms of the
canonical complex structure:

\begin{lemma} For a Frenet curve $f: M \to\HP^n$ the Willmore energy
  is given by
\begin{equation}
\label{eq:Willmore energy} W(f) = 2\int_M < A\wedge *A>\,.
\end{equation}
\end{lemma}
\begin{proof}
  Since $L$ is a Frenet curve the {mixed structure} $\einhalb(\nabla +
  *\nabla S)$ stabilizes $L$. But $\hat D = \einhalb(\nabla + *\nabla
  S)|_L$ has $\hat D_- =\nabla'_-|_L = A|_L$ and satisfies for
  $\alpha\in\Gamma(L\invers), \ \psi\in\Gamma(L)$ the product rule
\[ 
<D\alpha,\psi> + <\alpha, \hat D\psi> = \einhalb(d<\alpha,\psi> +
*d<\alpha,S\psi>)\,,
\]
where $D$ is the holomorphic structure on $L\invers$. This equation
and \eqref{eq:dbardual} imply $\hat D = \dbar - Q_{L\invers}^*$ and
hence $A|_L = - Q_{L\invers}^*$. Since $f$ is a Frenet curve, $A$ has
image in $L$ so that
\[
 < Q_{L\invers}\wedge *Q_{L\invers}> = < A|_L\wedge *A|_L> =<
A\wedge *A>\,.
\]

\end{proof}

\begin{definition}[see \cite{osculates}]
  A Frenet curve $f: M \to\HP^n$ is called \emph{Willmore} if $f$ is a
  critical point of the Willmore energy under compactly supported
  variations of $f$ by Frenet curves where we allow the conformal
  structure on $M$ to vary.
\end{definition}

\begin{definition}
  The \thedef{energy functional} of $S: M \to \Zz :=\{ S\in\End(V)
  \mid S^2 = - I\}$ is given by
\begin{equation}
\label{eq:energy map} E(S) =  \einhalb \int_M < \nabla S \wedge *\nabla S>
= 2 \int_M <Q\wedge *Q> +  <A\wedge*A>\,. 
\end{equation}
A map $S: M \to \Zz$ is called \thedef{harmonic} if it is a critical
point of the energy functional. \\

\end{definition}

Let $\nabla = \dbar + \partial + Q +A$ the decomposition of the
trivial connection $\nabla $ on $V$ with respect to a complex
structure $S$.  By changing the complex structure to $-S$ we get
$K\End_-(V) = K\Hom_+(\bar V, V)$ and $\partial$ and $\dbar$ on $V$
induce by \eqref{eq:dbar} a complex holomorphic structure on
$K\End_-(V)$. If we change the complex structure on $\bar K \End_-(V)$
to $-S$ then $\partial$ and $\dbar $ give a complex holomorphic
structure $\dbar$ on $K\End_-(\bar V) = \overline{\bar K \End_-(V)}$,
i.e., an antiholomorphic structure $\partial$ on $\bar K \End_-(V)$.

As in \cite[Prop.\ 5]{coimbra} one shows
\begin{theorem} 
\label{th:S harm}
Let $S: M \to \Zz$. Then following are equivalent
\begin{enumerate}
\item $S$ is harmonic.
\item $*Q$ is closed which due to \eqref{eq:dnabla*A} is the same as
  $*A$ is closed.
\item $Q$ is antiholomorphic, i.e., $\partial Q = 0$.
\item $A$ is holomorphic, i.e., $\dbar A = 0$.
\end{enumerate}
Moreover, if $f: M \to\HP^n$ is a Frenet curve and $S: M \to \Zz$ its
canonical complex structure, then $S$ is conformal, i.e.,
\[ 
<*\nabla S,*\nabla S> = <\nabla S, \nabla S>\,.
\]
\end{theorem}

The \emph{degree} of a complex quaternionic vector bundle $V= E \oplus
E$ is defined by the degree of the complex vector bundle $E$ which is
given by the $+i$--eigenspace of $S$, see \cite[Sec. \ 2.1]{Klassiker}.  Since
$\dbar + \partial$ is a complex connection on $V$, the degree of $V$
can be computed by
 \begin{equation}
\label{eq:degree}
2\pi \deg(V,S) = \int_M <S R^{\dbar + \partial}>  
        \stackrel{\eqref{eq:r2}}{=}  \int_M  <A \wedge *A>
- <Q\wedge *Q>\,.
\end{equation}
Combining \eqref{eq:Willmore energy}, \eqref{eq:energy map}, and
(\ref{eq:degree}),  we obtain:
\begin{cor} Let $f: M \to\HP^n$ be a Frenet curve with 
  canonical complex structure $S$. Then
\[
E(S)+ 4\pi \deg(V,S) = 2 W(L)\,.
\]
\end{cor}

Similar techniques as used in the $S^4$--case \cite[Thm.\ 3]{coimbra} give the
usual relation between the Willmore condition and harmonicity.

\begin{theorem}[see \cite{osculates}]
\label{t:harmonic}
A Frenet curve $f: M \to\HP^n$ is Willmore if and only if the
canonical complex structure of $f$ is harmonic, i.e.,
\[
d^\nabla*A=0\,.
\]
\end{theorem}

We have the Kodaira correspondence between holomorphic curves $f: M
\to\HP^n$ and base point free linear systems $H\subset H^0(L\invers)$.
For a Willmore curve $f: M\to\HP^n$, it is natural to ask for which
choices of basepoint free linear systems $\check H\subset
H^0(L\invers)$ the induced holomorphic curve $\check L\subset \check
H\invers$ is again Willmore.

\begin{prop}
\label{prop:linearsystem}
Let $f: M \to\HP^n$ be a Willmore curve. Let $L\subset V$ and
$H\subset H^0(L\invers)$ be the corresponding line bundle and
basepoint free linear system. Let $\check H\subset H^0(L\invers)$ be a
linear system with $H = V\invers \subset \check H$ so that the map
$\check f: M \to\HP^m$ given by the Kodaira correspondence has a
canonical complex structure which extends continuously into the
Weierstrass points. Then $\check f$ a Willmore curve in $\HP^m$ where
$ m = \dim \check H$.
\end{prop}
\begin{proof}
  Let $\check f_t: M \to\HP^m$ be a variation of $\check f$ so that
  the compact support $K$ does not contain Weierstrass points. Without
  loss of generality, we can assume that $\check f_t$ is unramified on
  $K$. Then $\pi: \check V = \check H\invers \to V$ defines a
  variation of $f$ by Frenet curves $f_t: M \to\HP^n$ by $\pi(\check
  L_t) = L_t$. Since the Willmore energy only depends on the
  holomorphic structure on $L\invers$ and not on the linear system, we
  see that
\[
\frac{\partial}{\partial t}W(\check f_t) =
\frac{\partial}{\partial t}W(f_t) =0\,.
\]
The usual arguments, see \cite{osculates}, show that the canonical
complex structure of $\check f$ is harmonic on $K$, i.e.,
\[
d^\nabla*A =0 
\]
away from the Weierstrass points. A recent result on the removability
of singularities of harmonic maps \cite{helein} shows that the
canonical complex structure extends smoothly into the Weierstrass
points, and therefore $\check f$ is a Frenet curve.
\end{proof}

\begin{example}  In \cite{paule} this construction is used to show
  that Willmore spheres $f: S^2 \to S^4$ are soliton spheres.  More
  precisely, there exists a 3--dimensional linear system $H\subset
  H^0(L\invers)$ such that the Kodaira embedding of $L$ into
  $H\invers$ is the dual curve of a twistor projection of a
  holomorphic curve in $\CP^{5}$.
\end{example}

In general, projections of Willmore curves $L\subset V$ into flat
subbundles $\check V\subset V$ fail to be Willmore, see \cite{paule}.

\begin{prop}
\label{p:projection} 
Let $f: M \to\HP^n$ be a Willmore curve with canonical complex
structure $S$ and let $H\subset H^0(L\invers)$ be the corresponding
linear system.  Let $\check H$ be an $S$ stable basepoint free linear
system $\check H\subset H \subset H^0(L\invers)$ with $m = \dim \check
H \ge 2$. Then
\[
\check L = \pi(L) \subset \check V = \check H\invers
\]
defines a Willmore curve $\check f: M \to\HP^m$.  Here $\pi: V =
H\invers \to \check V$ is the canonical projection.
\end{prop}

\begin{proof} Since $f$ is a holomorphic curve, the line bundle $L$ is
  full, i.e., $L$ is not contained in a lower dimensional flat
  subbundle of $V$. The kernel $\check H^\perp = \ker \pi$ of $\pi$ is
  $\nabla$ stable which shows that $\pi|_L \not = 0$.  Since $\check
  H$ is a linear system the induced connection $\check\nabla$ on
  $\check V$ satisfies $\pi \nabla = \check\nabla \pi$. Moreover,
\[
\pi S =: \check S \pi
\]
defines a complex structure on $\check V$ since $\ker \pi = \check
H^\perp$ is $S$ stable.  The complex holomorphic structures $\check
\nabla_+''$ and $\nabla_+''$ on $\check V$ and $V$ given by the
complex structures $\check S$ and $S$  are related by
\[
\check \nabla_+'' \pi = \pi \nabla_+''\,.
\]
Since $\nabla_+''$ and $\check\nabla_+''$ stabilize $L$ and $\pi L$
respectively, the map $\pi|_L $ is a complex holomorphic map. In
particular, the zeros of $\pi|_L$ are isolated and the complex bundle
$\Im \pi|_L$ can be extended smoothly across the zeros.  In other
words, $\Im \pi|_L$ defines a complex quaternionic line bundle $\check
L$. Note that $\check L_p= \pi L_p$ away from the isolated zeros of
$\pi|_L$.
 
  Let $L \subset V_1 \subset \ldots \subset V$ be the Frenet flag of
  $f$. Since $\pi \pi_L = \pi_{\check L} \pi$ we see 
\[
\check \delta_0
  \pi|_L = \pi \delta_0\,.
\]
If $\check \delta_0 = 0$ then $V_1$ is contained in the flat bundle $L
+ \ker \pi$ which has rank $\le n$ since $\dim \ker \pi = \rank V -
\rank \check V \le n -1$.  This contradicts the assumption that $L$ is
a full curve in $V$, i.e., the assumption that $\delta_k\not= 0$ for
$k=0,\ldots,n-1$.  Thus the  map
$\check \delta_0\not=0$ is complex holomorphic since 
\[
*\check\delta_0 = \check S \check \delta_0 = \check \delta_0 \check
S\,,
\]
and defines a vector bundle $\check V_1$.  Clearly, $\check V_1$
extends $\pi V_1$.
  
  Proceeding inductively, we see that $\check \delta_k \pi|_{V_k} = \pi
  \delta_k$ and $\check \delta_k \not = 0$ for all $0 \le k \le \rank
  \check V -2$. In particular, $\check L$ is a full curve in $\check
  V$ with Frenet flag $\check V_k = \pi V_k$. Moreover, $*\check
  \delta_k = \check S \check \delta_k = \check \delta_k \check S$
  yields that $\check S$ is an adapted complex structure.

  By construction $\check A = \einhalb *(\check \nabla \check S)'$ and
  $A = \einhalb *(\nabla S)'$ satisfy $\check A \pi = \pi A$, hence
  $\check S$ is the canonical complex structure of $\check f$.
  In particular $\check f$ is a Frenet curve, and 
\[
d^{\check\nabla}*\check A \pi = \pi d^\nabla *A =
  0\,.
\]
 shows that  $\check f$ is  Willmore.
\end{proof}

 \begin{rem}
\label{r:projection}
If $\dim\check H=1$ the same arguments as in the proof above show that
$(\pi(L), \pi S, \pi \nabla)$ defines a flat complex quaternionic line
bundle.
\end{rem}

Since the Hopf fields $A$ and $Q$ of a Willmore curve are holomorphic,
the zeros of $A$ and $Q$ are isolated. Therefore, Lemma
\ref{l:AaufLnull} implies:
\begin{cor} \label{cor:discrete sets}
  Let $S$ be the canonical complex structure of a Willmore curve $f:
  M\to\HP^n$.
\begin{enumerate}
\item If $A\not = 0$ then the set
\[
 \tilde M:= \{p\in M \mid L_p \subset \ker A_p\}
\]
has no inner points.
\item If $Q\not = 0$ then the set
\[
  \hat M:= \{p\in M \mid \Im Q_p \subset (V_{n-1})_p\}
\]
 has no inner points.
\end{enumerate}
\end{cor}

We collect some examples of Willmore curves in $\HP^n$ and  methods to
construct new Willmore curves out of given ones.
\begin{examples}
\label{ex:willmore}
\begin{enumerate}
\item Let $h: M\to\CP^{2n+1}$ be a (complex) holomorphic curve
  whose $n^{\text{th}}$ osculating space $W_n$ does not contain a
  quaternionic subspace, i.e, $W_n \oplus W_n j = \C^{2n+2} =
  \H^{n+1}$.  It is shown in \cite[Lemma 2.7]{Klassiker} that the
  \emph{twistor projection} $f: M \to\HP^n$ of $h$ has the smooth
  canonical complex structure $S$ given by $S|_W = i$. Moreover, it is
  shown that $A|_L=0$ so that Lemma \ref{l:AaufLnull} gives that $f$
  is Willmore since $A=0$. Moreover, $f$ has Willmore energy $W(f)=0$.
  
  Conversely, every Willmore surface $f: M \to\HP^n$ with $W(f)=0$ is
  given as a twistor projection of a holomorphic curve in
  $\CP^{2n+1}$.
  
\item Let $f: M \to\HP^n$ be a Willmore curve, and $L\subset V$ the
  corresponding line bundle. The flat connection $\nabla^*$ on
  $V\invers$ decomposes as $\nabla^* = \dbar^* + \partial^* - A^* -
  Q^*$ with respect to the canonical complex structure $S^*$ of the
  dual curve $f^*: M \to\HP^n$. Therefore, we compute
\[
(d^{\nabla^*} *A^*)^* = d^\nabla *A
\] 
so that a curve $f: M\to\HP^n$ is Willmore if and only if its dual
curve $f^*: M \to\HP^n$ is Willmore.
  
\item A flat connection $\nabla$ on a complex quaternionic vector
  bundle $(V,S)$ is called \thedef{Willmore connection} \cite[Sec.\ 
  6.1]{Klassiker} if $S$ is harmonic, i.e., $d^\nabla*A = 0$. In
  general, the harmonic complex structure $S$ will not be the
  canonical complex structure of a Frenet curve.  But if $\rank A = 1$
  then $\dbar A = 0$ implies that the image of $A$ defines a $\dbar$
  holomorphic line bundle.  Thus, $L = \Im A \subset V$ is a Willmore
  curve if $Q\delta^k = 0$ for all $k=0,\ldots, n-1$ see Lemma
  \ref{lem:QA=0}.
  
  Moreover, if $L$ is Willmore then the connections $\nabla^\lambda =
  \nabla + (\lambda -1)A$ are flat for all $\lambda = \alpha + \beta
  S, \ \alpha,\beta\in \R, \ \alpha^2 + \beta^2 =1.$ Denote by
  $L^\lambda$ the line bundle $L$ considered as subbundle in
  $(V,\nabla^\lambda)$. If we decompose $\nabla^\lambda$ with respect
  to the complex structure $S$ of $L$ then $Q^\lambda = Q$. Thus $S$
  is the canonical complex structure of $L^\lambda$ and $L^\lambda$ is
  a Willmore curve. Its Willmore energy is given by $W(L^\lambda) =
  W(L)$.
  
  Notice, that though $\nabla$ is trivial, the Willmore curves of this
  family may have holonomy.
\end{enumerate}
\end{examples}



\section{Pl\"ucker relation of a Frenet curve}

The quaternionic Pl\"ucker relation \cite[Thm.\ 4.7]{Klassiker} gives the
Willmore energy of a Frenet curve in terms of the Willmore energy of
its dual curve, the genus of the surface and the degree of the
associated line bundle. We give a proof of the Pl\"ucker relation in
the case when $f: M \to\HP^n$ is a Frenet curve.

We compute the degrees of various complex bundles involved in the
Pl\"ucker relation:
\begin{lemma}
\label{lemma:dualdegree}
Let $f: M \to\HP^n$ be a Frenet curve, $S$ its canonical complex
structure and $L\subset V_1 \subset \ldots \subset V_n = V$ its Frenet
flag with corresponding derivatives $\delta_i$. Then the degree of the
bundle $V_k/V_{k-1}$ with respect to $S$ is given by
\begin{equation}
\label{eq:dualdegree}
\deg V_k/V_{k-1}  = \sum_{i=0}^{k-1}\ord \delta_i - k \deg K + \deg L,
\ 0 \le k \le n\,.
\end{equation}
where we put $V_{-1}:=\{0\}$.
\end{lemma}
\begin{proof}
  The degree of a complex holomorphic line bundle $E$ is given by the
  vanishing order of any holomorphic section of $E$. Since $*\delta_i
  = S\delta_i = \delta_i S$, Lemma \ref{lem:holbundel} implies that
  $\delta_i\in H^0(K\Hom_+(V_i/V_{i-1}, V_{i+1}/V_i))$ is
  a  holomorphic section, and thus
\[
\ord\delta_i = \deg(K\Hom_+(V_i/V_{i-1}, V_{i+1}/V_i))\,.
\]
If $(V_1,S_1), (V_2,S_2)$ are two complex quaternionic vector bundles
then $\Hom_+(V_1,V_2)$ is canonically isomorphic to $\Hom_\C(E_1,
E_2)$, where the $E_k$ are again the $+i$--eigenspaces of $S_k$.
Therefore, 
\[
\ord \delta_i = \deg K + \deg V_{i+1}/V_i - \deg V_i/V_{i-1}\,,
\]
and, 
telescoping this identity, we get
\begin{equation}
\label{eq:degreeofquotients}
 \sum_{i=0}^{k-1}\ord \delta_i =  k \deg K + \deg V_k/V_{k-1}  - \deg L\,.
\end{equation}
\end{proof}

\begin{rem} The degree of the dual curve is given by
\[ 
\deg L^* = n\deg K - \deg L -\sum_{i=0}^{n-1}\ord \delta_i\,,
\]
since $(V/V_{n-1})\invers = V_{n-1}^\perp = L^*$.
\end{rem}

\quad

We now prove the quaternionic Pl\"ucker relation \cite[Thm.\ 
4.7]{Klassiker} in the case of Frenet curves:
\begin{theorem}
\label{thm:Pluecker}
Let $f: M \to\HP^n$ be a Frenet curve with canonical complex structure
$S$. Let $L \subset V_1 \subset \ldots \subset V_n = V$ be the Frenet
flag of $f$ and $\delta_i$ the derivatives of $V_i$.  For a compact
Riemann surface $M$ of genus $g$, \emph{the Pl\"ucker relation}
\begin{equation}
\label{eq:Pluecker}
\deg(V,S) =\frac{1}{4\pi}(W(f) - W(f^{*})) = (n+1)(n(1-g) + \deg L)
  +\ord
 H
\end{equation}
holds, where $\ord H= \sum_{i=0}^{n-1} (n-i) \ord \delta_i$ is the order of the linear system $H =
V\invers\subset H^0(L\invers)$.
\end{theorem}
\begin{rem} If $L\subset V$ is a  holomorphic curve then $H=
  V\invers\subset H^0(L\invers)$ is a basepoint free linear system.
  The \emph{order} of $H$ is defined by \cite[Def.\ 4.2]{Klassiker}
\[ 
\ord(H)=\sum_{p\in M}\ord_p(H)\,.
\]
where $\ord_p(H) = \sum_{k=0}^n(n_k(p) - k)$ is the order of $H$ at
$p$ and $n_0(p) < \ldots < n_n(p)$ is the Weierstra\ss \ gap sequence
of $H$.  In the case of a Frenet curve, the expression for the order
of $H$ simplifies to $\ord H= \sum_{i=0}^{n-1} (n-i) \ord \delta_i$.
In particular, if $f$ is an unramified Frenet curve, then $\ord H=0$.
\end{rem}
\begin{proof}
  Since as complex vector bundles $V = \bigoplus_{k=0}^n V_k/V_{k-1}$
  we have
\begin{eqnarray*}
\deg(V,S) &=& \sum_{k=0}^n \deg V_k/V_{k-1}\stackrel{\eqref{eq:degreeofquotients}}{=}  \sum_{k=0}^n(
\sum_{i=0}^{k-1}\ord \delta_i -  k \deg K + \deg L) \\
&=& (\sum_{k=0}^n
\sum_{i=0}^{k-1}\ord \delta_i) -  \frac{n(n+1)}{2}\deg K + (n+1)\deg
L\,.
\end{eqnarray*}
Moreover,
\[ 
4\pi \deg(V,S)  \stackrel{\eqref{eq:r2}}{=}
 2 \int_M  <A \wedge *A> - <Q\wedge *Q> 
\stackrel{\eqref{eq:Willmore energy}, \   \eqref{eq:Adagger}}=
 W(f) - W(f^*)\,.
\]
\end{proof}

\begin{rem}
\label{r:Willmoreenergyoftwistor}
Let $h: M \to\CP^{2n+1}$ be a holomorphic curve in $\CP^n$ such that
the twistor projection $f: M\to\HP^n$ of $h$ is a Frenet curve,
compare Examples \ref{ex:willmore} \thetag{1}. Since $W(f)=0$ the
Pl\"ucker relation shows that the Willmore energy of the dual curve
$f^*$ of $f$ is given by
\[
 W(f^*)  = 4\pi \deg(V,S) \in 4\pi\N\,.
\]
\end{rem}



\section{B\"acklund transformation on Willmore curves }

Using the harmonicity of the canonical complex structure, a similar
construction to the $\delbar$ and $\del$ transforms of a harmonic map
into $\CP^n$, \cite{wolfson}, gives the B\"acklund transformation on
Willmore curves. We show that the B\"acklund transform $\tilde f:M
\to\HP^n$ of a Willmore curve $f: M \to\HP^n$ is again Willmore
provided $\tilde f$ is a Frenet curve. The latter assumption will be
void in case of Willmore spheres.

Due to the harmonicity of the canonical complex structure
$A$ is $\dbar$--holomorphic and $Q$ is $\partial$--holomorphic.  Thus
their kernels and images define smooth subbundles of the trivial
$\H^{n+1}$--bundle, and we get new maps into $\HP^n$. 

\begin{lemma}
\label{l:hatL}
Let $f: M \to\HP^n$ be a Willmore curve and $S$ its canonical complex
structure. Decomposing the trivial connection $\nabla = \dbar +
\partial + Q + A$ on $V=\trivial{n+1}$ with respect to $S$, we see
\begin{enumerate}
\item For $ A \not = 0 $ there exists a rank $n$ subbundle $\tilde
  W_{n-1}\subset V$ which agrees with $\ker A$ except at finitely many
  points and which satisfies $\tilde W_{n-1} \subset \ker A$.
  
\item For $Q\not = 0$ there exists a line bundle $\hat L\subset V$
  which agrees with $\Im Q$ except at finitely many points and
  satisfies $\Im Q\subset \hat L$.
\end{enumerate}  
In abuse of notation, we write $\ker A = \tilde W_{n-1}$ and $\Im Q =
\hat L$.
\end{lemma}
The image of $A$ and the kernel of $Q$ also define smooth bundles.
However, these are the already known bundles $\Im A= L$ and $\ker Q
=V_{n-1}$ since $S$ is the canonical complex structure of $f$.  The
harmonicity of $S$ implies that $\ker A$ is a holomorphic bundle and
$\Im Q$ a holomorphic curve.
\begin{lemma}
\label{l:tildeL}
Let $f: M \to\HP^n$ be a Willmore curve with canonical complex
structure $S$.
\begin{enumerate} 
\item If $A\not = 0$ then $\ker A \subset V$ is a holomorphic
  subbundle with respect to the holomorphic structure induced by the
  complex structure $-S$ on $V$.
\item If $ Q\not= 0$ then $\Im Q\subset V$ is a holomorphic curve with
  respect to the complex structure $-S$ on $V$.
\end{enumerate}
\end{lemma}
\begin{proof}
  For $\varphi \in \Gamma(\ker A)$ we have
 \begin{eqnarray*}
  0 &=&(d^\nabla *A)\varphi = d^\nabla(*A\varphi) + *A\wedge\nabla\varphi  
     = *A*\nabla\varphi + A\nabla \varphi \\ 
    &=& *A*\tilde\delta_{n-1} \varphi + A\tilde\delta_{n-1}\varphi
     = -AS(*\tilde\delta_{n-1} \varphi +S\tilde\delta_{n-1}\varphi)\,.  
 \end{eqnarray*}
which implies 
\begin{equation}
\label{eq:kerAholomorph}
*\tilde\delta_{n-1}  +
S\tilde\delta_{n-1} = 0\,,
\end{equation} 
since $A$ is, interpreted as a map in $\Omega^1(\Hom(V/\ker A,V))$,
injective away from finitely many points.  This shows that $\ker
A\subset V$ is a holomorphic subbundle, compare Example
\ref{ex:holo_structure}\thetag{3}.

Consider the dual curve $f^*$ of $f$ which is again a Willmore curve
with $A^\dagger = - Q^*$. Assume that $Q \not = 0$, then by the above
argument $\ker A^\dagger\subset V\invers$ is a holomorphic subbundle.
But then $\Im Q\subset V$ is a holomorphic curve, by Example
\ref{ex:holo_structure}\thetag{3} and
\begin{equation}
  \label{eq:dual}
  (\ker A^\dagger)^\perp =(\ker Q^*)^\perp = \Im Q\,.
\end{equation}
\end{proof}

To be able to deal with holomorphic curves only, we consider instead
of the holomorphic bundle $\ker A$ the holomorphic curve $(\ker
A)^\perp\subset V\invers$.  One of the main difficulties of this
construction is that $\Im Q$ and $(\ker A)^\perp$ might fail to be
Frenet curves.  We will show below that at least in the case of
Willmore spheres $f: S^2 \to\HP^n$ both line bundles are Frenet
curves. The general case is more difficult and is a topic to which we
will return in a future paper.

Therefore, at least for the purposes of the present paper, we will
assume that $\hat L$ and $(\ker A)^\perp$ are Frenet curves. We define
\[
\hat V\subset V \quad \text{and} \quad \tilde V\invers\subset V\invers
\]
as the trivial subbundles of $V$ and $V\invers$ so that the
holomorphic curves $\hat L\subset \hat V$ and $(\ker A)^\perp\subset
\tilde V\invers$ are full curves in $\hat V$ and $\tilde V\invers$
respectively.

Since $(\ker A)^\perp$ is a line subbundle of the dual bundle of $V$,
we will rather consider the dual curve $\tilde L\subset \tilde V$ of
$(\ker A)^\perp$ which is again a Frenet curve unless $(\ker A)^\perp$
is a constant in $\P V$. In this case, we define $\tilde L:= ((\ker
A)^\perp)\invers$ to be the dual bundle of $(\ker A)^\perp$.

\begin{definition} Let $f: M \to\HP^n$ be a Willmore curve. Then 
  $\tilde L\subset \tilde V$\, is called the  \emph{forward
    B\"acklund transform} of $L$, and $\hat L\subset \hat V$ the
  \emph{backward B\"acklund transform} of $L$.
\end{definition}

\begin{rem}
  Note the similarity to the $\delbar$ and $\partial$ transforms of
  harmonic maps into $\CP^n$ \cite{wolfson}: we use the $(0,1)$--part
  $Q$ and the $(1,0)$--part $A$ of the derivative $\nabla S$ of the
  harmonic map $S: M\to\mathcal Z$ to construct new holomorphic
  curves.  However, our construction will give new Willmore curves
  rather than the associated harmonic maps.  Moreover, to obtain
  sequences of Willmore surfaces, we will have to guarantee the
  smoothness of the canonical complex structure of a B\"acklund
  transform. This can be done, at least in the case of Willmore
  spheres in $\HP^n$, and we will see below that in this case the
  resulting sequence is finite.
  
  To compare our definition of the B\"acklund transformation to the
  one given in \cite[Prop.\  17]{coimbra} for conformal immersions $f: M
  \to S^4$, we contemplate the Frenet flag of $\tilde L$.
  
  Let $k=\rank \tilde V$ be the rank of the trivial bundle $\tilde V$
  and let $\tilde V_i^*$ be the Frenet flag of the Frenet curve $
  (\ker A)^\perp \subset \tilde V\invers$. The Frenet flag of the
  B\"acklund transform $\tilde L = ((\ker A)^\perp)^*$ is thus given
  (\ref{eq:dualflag}) by
\[
\tilde V_i = (\tilde V^*_{k-1-i})^\perp\,.
\]
In the case when $\tilde L$ is a full curve in $V$ then
\[
\tilde V_{n-1} = (\tilde V^*_{0})^\perp = (\ker A)^{\perp\perp} =
\ker A\,.
\]
In particular, for Willmore surfaces $f: M\to S^4$ the B\"acklund
transform $\tilde f$ is a full curve in $\H^2$ unless it is a constant
point in $\HP^1$. Therefore, we obtain the (twofold) forward
B\"acklund transform $\tilde L = \tilde V_{n-1} =\ker A$ as defined in
\cite[Prop.\  17]{coimbra}.
\end{rem}

We prove that B\"acklund transforms of Willmore curves are again
Willmore curves:
\begin{theorem} The forward and the backward B\"acklund transform of
  a Willmore curve are again Willmore curves.
\end{theorem}
\begin{proof}
  Let $f: M \to\HP^n$ be a Willmore curve, $S$ its canonical complex
  structure and $\tilde f: M \to\HP^k$ the forward B\"acklund
  transform of $f$. The line bundle $(\ker A)^\perp =: \tilde L^*$ is
  a Frenet curve in some trivial quaternionic subbundle $\tilde
  V\invers\subset V\invers$ of $\rank k+1$. We denote the induced
  projection by $\pi: V \to \tilde V$.  The trivial connection
  $\nabla$ on $\tilde V\invers$ induces a trivial connection $\nabla$
  on $\tilde V$. Since $(\tilde L^*)^\perp$ equals $ \pi \ker A$, the
  Frenet flag of $\tilde f$ is given by $\tilde L \subset \tilde V_1
  \subset \ldots \subset \tilde V_{k-1} \subset \tilde V$ where
  $\tilde V_{k-1} = \pi \ker A$.
  
  Since $\ker A$ is a holomorphic vector bundle with respect to the
  holomorphic structure induced by $-S$, the line bundle $\tilde L^*$
  is a holomorphic curve with respect to the complex structure
  $-S^*|_{\tilde V\invers}$.  In particular, the canonical complex
  structure of $\tilde f^*$ is given by
\[
\tilde S^* =-S^*|_{\tilde
    V\invers} + \tilde B^*
\]
with $\tilde L^* \subset \ker \tilde B^* = (\Im \tilde B)^\perp$.
Therefore
\[ 
\tilde S \pi = - \pi S + \tilde B\,,
\]
defines the canonical complex structure of $\tilde f$ where $\tilde
B\in\Gamma(\Hom(V, \tilde V_{k-1}))$.  The bundle $\tilde V_{k-1}$ is
$\tilde S$ stable and for $\varphi\in\Gamma(V)$ we calculate
\begin{eqnarray*}
\pi_{\tilde V_{k-1}}(*\nabla \tilde B - \tilde S\nabla \tilde
 B)\varphi
 &=& 
\pi_{\tilde V_{k-1}}(*\nabla (\tilde B\varphi) -
\tilde B*\nabla \varphi - \tilde S\nabla(\tilde B\varphi) + \tilde S \tilde B \nabla \varphi)\\
&=& 
\pi_{\tilde V_{k-1}}(*\nabla (\tilde B\varphi) - \tilde S\nabla (\tilde B\varphi
)) = (*\tilde\delta_{k-1} - \tilde S\tilde\delta_{k-1})\tilde B\varphi
 \\
 &=& 0\,.
\end{eqnarray*}
This shows that $*\nabla \tilde B - \tilde S\nabla \tilde B \ \text{
  takes values in } \ \tilde V_{k-1}$. Since $Q|_{V_{n-1}}=0$ we also
obtain
\begin{eqnarray*}
4 \tilde A \pi |_{V_{n-1}} &\stackrel{\eqref{eq:*A}}=& 
(\tilde S\nabla\tilde S + *\nabla\tilde S) \pi |_{V_{n-1}}\\
&=&  
\big(\pi\big (S(\nabla S) - *\nabla S\big ) - \tilde B(\nabla S) + \tilde S
(\nabla \tilde B) + *\nabla \tilde B\big)|_{V_{n-1}} \\
&\stackrel{\eqref{eq:*A}}=& 
 \big( 4\pi Q - \tilde B(\nabla S) + \tilde S
(\nabla \tilde B) + *\nabla \tilde B\big)|_{V_{n-1}} \\
&=&
\big (-\tilde B(\nabla S) + \tilde S
(\nabla \tilde B) + *\nabla \tilde B\big)|_{V_{n-1}}\,. 
\end{eqnarray*}

But $\tilde A$ maps to $\tilde L \subset \tilde V_{k-1}$, so we see
that $*\nabla \tilde B + \tilde S \nabla \tilde B$ restricted to
$V_{n-1}$ takes values in $\tilde V_{k-1}$ and so does $\nabla \tilde
B$. For $\psi \in V_{n-1}$ we get
 \[
\tilde\delta_{k-1}\tilde B\psi = \pi_{\tilde V_{k-1}}(\nabla\tilde B)\psi = 0
\]
 and hence
\[ 
\tilde S \pi = - \pi S + \tilde B, \text{ where } \Im \tilde B\subset
\tilde V_{k-1} \text{ and } \ V_{n-1} \subset \ker \tilde B\,.
\] 
This yields 
 \begin{equation*}
\label{s5} 
*\nabla \pi \psi + (\nabla\tilde S)\pi\psi +  \tilde S \nabla \pi \psi 
=
  \pi *\nabla \psi +\nabla (\tilde S\pi \psi)
 = 
\pi(*\nabla \psi - \nabla S\psi) \in\Omega^1(\pi L)\,,
\end{equation*}
for $\psi \in \Gamma(V_{n-1})$ since $\pi_L(*\nabla \psi - \nabla S
\psi) = *\delta_0\psi - \delta_0 S\psi = 0$. Moreover, $\tilde S$
stabilizes $\pi L$ by $\tilde S\pi \varphi = - \pi S\varphi + B\varphi
= - \pi S\varphi \in \pi L$ for $\varphi\in L$.  Thus we also get
\begin{equation*}
\tilde S *\nabla \pi \psi + \tilde S(\nabla \tilde S)\pi\psi - \nabla \pi\psi
\in\Omega^1(\pi L)
\end{equation*}
and
\begin{equation*}
-\nabla \pi\psi + (*\nabla \tilde S)\pi\psi + \tilde S *\nabla \pi\psi
\in\Omega^1(\pi L)\,.
\end{equation*}
Subtracting these equations, we find
\[ 
4\tilde Q \pi\psi \stackrel{\eqref{eq:*A}}= 
(\tilde S(\nabla \tilde S)- *\nabla \tilde S)\pi \psi\in\Omega^1(\pi
L) 
\ \text{ for } \psi \in \Gamma(V_{n-1})\,.
\]
Since $\tilde S$ is the canonical complex structure of $\tilde f$,
i.e., \[
\tilde Q|_{\pi(\ker A)}= \tilde Q|_{\tilde V_{k-1}}=0\,,
\]
this implies that $\tilde Q $ takes values in $\pi L$.  Moreover,
\begin{eqnarray*}
4\tilde Q \pi &\stackrel{\eqref{eq:*A}}=& (\tilde S \nabla \tilde S - *\nabla \tilde S)\pi \\
&=& (\tilde S\nabla \tilde B - \tilde B\nabla S - *\nabla \tilde B) + \pi(S\nabla S + *\nabla S)\\
&\stackrel{\eqref{eq:*A}}=&(\tilde S\nabla \tilde B - \tilde B\nabla S - *\nabla \tilde B)
+  4 \pi A\,.
\end{eqnarray*}
Since $\tilde B$ and $ *\nabla \tilde B - \tilde S\nabla \tilde B$ map
to $\tilde V_{k-1}$ while $\pi A, \tilde Q$ have values in $ \pi L$,
we obtain $(\tilde S\nabla \tilde B - \tilde B\nabla S - *\nabla
\tilde B) = 0$ and 
\begin{equation}
\label{eq:fQ} \tilde Q \pi  = \pi A\,.
\end{equation}
Now $(d^\nabla *\tilde Q)\pi = \pi (d^\nabla *A)=0$ yields by Theorem
\ref{t:harmonic} that $\tilde f$ is Willmore. \\

Assume that $f: M \to\HP^n$ is a Willmore curve such that $Q\not=0$
and such that $\hat f$ is a Frenet curve in some trivial quaternionic
subbundle $\hat V \subset V$. The dual curve $f^*$ of $f$ is Willmore,
 and by the above argument $\widetilde{f^*}$ is Willmore, too. Finally,
\begin{equation}
\label{eq:tilde/hat}
\widetilde{L^*} = ((\ker A^\dagger)^\perp)^* \stackrel{\eqref{eq:dual}}=
 \hat L^*\,.
\end{equation}
shows that $\hat f$ is the dual curve of $\widetilde{f^*}$ and
therefore Willmore. 

For later use we collect the information we have on the canonical
complex structure $\hat S$ of $\hat f$: The canonical complex structure
$\widetilde{S^*}$ of $\widetilde{f^*}$ induces $\hat S$ via
$
 \hat S = \widetilde{S^*}^*\,,
$ and thus 
\[
\hat S = -S|_{\hat V} + \hat B\,,
\] 
where $\hat L \subset \ker \hat B$ and $\Im \hat B \subset L$.  Let
$\pi: V\invers \to \hat V\invers$ be the canonical projection. Then
$\pi A^\dagger = \widetilde{Q^\dagger} \pi $ and
\begin{equation}
\label{eq:bA} \hat A = - (\widetilde{Q^\dagger}\pi )^* = - (\pi
A^\dagger)^* =  Q|_{\hat V}\,.
\end{equation}
\end{proof}

As a consequence of  \eqref{eq:tilde/hat}, we get the relation
between the forward B\"acklund transform of a Willmore curve and the
backward B\"acklund transform of its dual curve. The dual statement
follows similarly.
 
\begin{cor}\label{cor:fdual=dualb}
  Let $f: M \to\HP^n$ be a Willmore curve and $f^*$ its dual curve.
\begin{enumerate}
\item If the forward B\"acklund transform of  $f^*$ exists, then the
  backward B\"acklund transform of $f$ exists and
\[ \widetilde{f^*} = \hat f^*.\]
\item If the backward B\"acklund transform of $f^*$ exists, then the
  forward B\"acklund transform of $f$ exists and
\[ \widehat{f^*}=\tilde f^*.\]
\end{enumerate}
\end{cor} 

\quad

If the B\"acklund transform $\tilde f$ of a Willmore curve $f$ is a
Frenet curve in $\HP^n$ then $\tilde f$ has $\tilde Q = A$ by
\eqref{eq:fQ}.  Hence $\Im \tilde Q = \Im A = L$ and the backward
transform of $\tilde f$ exists.

\begin{cor} Let $f: M\to\HP^n$ be a Willmore curve.
\begin{enumerate}
\item Assume that the forward B\"acklund transform $\tilde f$ is a
  Frenet curve in $\HP^n$. Then the backward B\"acklund transform of
  $\tilde f$ exists and
  \begin{equation*}
      \widehat{\tilde f} = f 
  \end{equation*}
\item Assume that the backward B\"acklund transform $\hat f$ is a
  Frenet curve in $\HP^n$. Then the forward B\"acklund transform
  of $\hat f$ exists and
  \begin{equation*}
       \widetilde{\hat f} = f.
  \end{equation*}
\end{enumerate}
\end{cor}

We are now able to give the Willmore energy of a B\"acklund
transform of $f$ in terms of the Willmore energy of $f$.

\begin{cor}
\label{cor:willmoreofbaecklund}
Let $f: M \to\HP^n$ be a Willmore curve, and let $\hat f$ and $\tilde f$
be the backward and forward B\"acklund transforms of $f$. Then the
Willmore energies of $\hat f$ and $\tilde f$ are given by
\[
 W(\hat f) = W(f^*) \quad \text{ and } \quad W(\tilde f^*) = W(f)\,.
\]
\end{cor}
\begin{proof}
\quad

\begin{enumerate}
\item Recall that by \eqref{eq:Adagger} and \eqref{eq:Willmore energy}
  the Willmore energy of the dual curve $f^*$ is given by
  $W(f^*)=\int_M <Q\wedge*Q>$. By \eqref{eq:bA} we have $\hat A =
  Q|_{\hat V}$. Since $\Im Q = \hat L \subset \hat V$ we get
\begin{eqnarray*}
\qquad  W(f^*)&=&\int_M <Q \wedge *Q> =\int_M < Q|_{\hat V}\wedge *Q|_{\hat
   V}> =\int_M<\hat A \wedge *\hat
A> = W(\hat f).
\end{eqnarray*}

\item Using Corollary \ref{cor:fdual=dualb} we see
$W(\tilde f^*)  = W(\widehat{f^*}) \stackrel{\text{\thetag{i}}}{=} W(f)$.
\end{enumerate}
\end{proof}



\section{B\"acklund transforms with $-S$ as the canonical complex structure}

Given the forward and backward B\"acklund transforms $\tilde f$ and
$\hat f$ of a Willmore curve $f: M \to\HP^n$, we have seen that the
negative $-S$ of the canonical complex structure of $f$ renders
$\tilde f^*$ and $\hat f$ into holomorphic curves.  We will now
discuss the case when $-S$ is in fact the canonical complex structure
of $\tilde f$ or $\hat f$. It turns out that in this case the
B\"acklund transform comes from complex holomorphic data and $f$ has
integer Willmore energy.

Moreover, the B\"acklund transforms can be used to project $f$ to a
Willmore curve $\check f: M \to\HP^{n-k}$ for some suitable
$\HP^{n-k}\subset \HP^n$ such that $\check f$ is given by complex
holomorphic data.

\begin{theorem}
  Let $f: M \to\HP^n$ be a Willmore curve.
\begin{enumerate}
\item If the backward B\"acklund transform $\hat f: M\to\HP^k$ is a
  Frenet curve in $\HP^k$ with $k\le n$ and has canonical complex
  structure $\hat S = -S$, then $\hat f$ is the dual curve of a
  twistor projection of a holomorphic curve $h: M \to\CP^{2k+1}$.
\item If the forward B\"acklund transform $\tilde f: M\to\HP^k$ is a
  Frenet curve in $\HP^k$ with $k\le n$ and has canonical complex
  structure $\tilde S = - S$ then $\tilde f$ is the twistor projection
  of a holomorphic curve $h: M \to\CP^{2k+1}$.
\end{enumerate}
In both cases, $f$ has Willmore energy $ W(f) \in 4 \pi \N$.
\end{theorem}

\begin{proof}
 
\quad

\begin{enumerate}
\item Let $\hat V$ be the trivial $k+1$ bundle so that $\hat L \subset
  \hat V$ is a full curve. Since $\hat S = - S|_{\hat V}$ is the
  canonical complex structure of $\hat f$, Lemma \ref{lem:QA=0} shows
  that
\[
A\hat\delta^k=0\,,
\]
which implies $A|_{\hat V}=0$. Since $\hat Q = A|_{\hat V}$ this
yields that $\hat f$ is the dual curve of the twistor projection of a
holomorphic curve in $\CP^{2k+1}$. By Remark
\ref{r:Willmoreenergyoftwistor} $\hat f$ has integer Willmore energy,
and the Pl\"ucker relation together with Corollary
\ref{cor:willmoreofbaecklund} gives
\[ 
W(f) = W(\hat f) + 4\pi\deg(V,S) \in 4\pi \N\,.
\]

\item The dual curve $f^*: M \to\HP^n$ of $f$ is a Willmore curve.  By
  Corollary \ref{cor:fdual=dualb} the backward B\"acklund transform of
  $f^*$ is given by the dual
\[
\widehat{f^*} = \tilde f^*
\]
of the forward B\"acklund transform of $f$. Since $\tilde f: M
\to\HP^n$ has canonical complex structure $-S$, the backward
B\"acklund transform $\widehat{f^*}$ has canonical complex structure
$\widehat{S^*} = - S^*$. Using the first part, $\tilde f =
(\widehat{f^*})^*$ is the twistor projection of a holomorphic curve
$h: M \to\CP^{2k+1}$, and $W(f) = W(\tilde f^*)\in 4\pi\N$.
\end{enumerate}
\end{proof}

For Willmore curves $f: M \to\HP^n$ the backward B\"acklund transform
$\hat f: M \to\HP^k$ is a Frenet curve in some trivial $\rank k+1$
subbundle $\hat V\subset V$ of $V$, $k\le n$. If the canonical complex
structure $\hat S$ is given by $-S$ on $\hat V$, then $\hat V$ is in
particular $S$ stable, too.

If $\rank \hat V = n+1$, i.e., $V = \hat V$, then $A\delta^i=0$ for
$i\le n$ since $-S$ is the canonical complex structure of $\hat f$.
Because $\Im \delta^i = V_{i+1}/V_i$ except at finitely many points,
this shows that $A =0$. In other words, if the backward B\"acklund
transform $\hat f$ is a full curve in $V$, then $f$ is the twistor
projection of a holomorphic curve in $\CP^{2n+1}$.

If $\rank \hat V = k +1 < n$, the quotient bundle $\check V = V/\hat
V$ is a smooth trivial bundle of $\rank n-k > 1$.  The canonical
projection $\pi: V\to \check V$ has $S$ stable kernel $\ker\pi = \hat
V$ so that we can project $L$ to a Willmore curve $\check L =
\pi(L)\subset \check V$ by Proposition \ref{p:projection}.  In other
words, $f: M\to\HP^n$ projects to a Willmore curve $\check f: M \to
\HP^{n-k}$.

Moreover, we know that the canonical complex structure of $\check f$
is given by $\check S \pi = \pi S$ so that
\[
\check Q \pi = \pi Q\,.
\]
The image $\hat L$ of $Q$ is a line subbundle of $\hat V$ and thus
$\pi Q =0$ and $\check Q=0$.  Therefore, we have shown that $\check f$
is the dual curve of a twistor projection of a holomorphic curve in
$\CP^n$. Dualizing this result, we obtain a similar result in case
that the forward B\"acklund transform has $-S$ as complex structure.

\begin{prop}
\label{prop:sphere projection}
Let $f: M \to\HP^n$ be a Willmore curve so that $f$ and $f^*$ do not
come from the twistor projection of a holomorphic curve $h: M
\to\CP^{2n+1}$. 
\begin{enumerate}
\item If the backward B\"acklund transform $\hat f: M \to\HP^k$, $k<
  n$, has canonical complex structure $-S$ then the line bundle $L$
  projects under the canonical projection $\pi: V \to\check V$ to a
  Willmore curve $\check L = \pi(L)$ in the rank $n-k$ trivial bundle
  $\check V = V/\hat V$ unless $\check L =\check V$ is
  1--dimensional. 

  More precisely, $\pi(f)= \check f: M \to\HP^{n-k}$ is the dual curve
  of a twistor projection of a holomorphic curve $h: M \to \CP^n$
  unless $\check f$ is a constant point. \\

\item If the forward B\"acklund transform $\hat f: M \to\HP^k$, $k<
  n$, has canonical complex structure $-S$ then the dual line bundle
  $L^*$ of $L$ projects under the canonical projection $\pi: V\invers
  \to\check{(V^*)}$ to a Willmore curve $\check{(L*)} = \pi(L^*)$ in
  the rank $n-k$ trivial bundle $\check{(V^*)} = V\invers/\tilde
  V\invers$ unless $\check{(L^*)} =\check{(V^*)}$ is
  1--dimensional. 

  More precisely, $\pi(f^*)= \check{(f^*)}: M \to\HP^{n-k}$ is the
  dual curve of a twistor projection of a holomorphic curve $h: M \to
  \CP^n$ unless $\check{(f^*)}$ is a constant point.
\end{enumerate}
\end{prop}

\begin{rem}
  The assumption that the canonical complex structure of a B\"acklund
  transform of a Willmore curve $f: M \to\HP^n$ is the negative of
  the canonical complex structure of $f$ does not hold in general: for
  example, there are Willmore tori $f: T^2\to S^4$ which come from
  integrable system methods and do not have integer Willmore energy,
  c.f. \cite{s4paper}.
\end{rem}



\section{Willmore spheres in $\HP^n$}

We show that the forward and backward B\"acklund transforms of a
Willmore sphere $f: S^2\to\HP^n$ are Frenet curves whose canonical
complex structures are the negative $-S$ of the canonical complex
structure $S$ of $f$. In particular, combining the results of the
previous section, we see that a Willmore sphere $f$ has integer
Willmore energy $W(f)\in 4\pi\N$ and is either a minimal surface in
$\R^4$ with planar ends, or $f$ or its dual curve $f^*$ is, at most
after projection to a suitable $\HP^m$, a twistor projection of a
holomorphic curve in complex projective space. In particular, a
Willmore sphere is given by complex holomorphic data.

To show that a B\"acklund transform of a Willmore sphere is Frenet, we
construct the Frenet flag and an adapted complex structure
recursively.

\begin{lemma} \label{l:flag}
  Let $f: M \to\HP^n$ be a Willmore curve with canonical complex
  structure $S$ such that $Q\not= 0$ and let $\hat L = \Im Q$.  Assume
  that there exists for $0\le k\le n-1$ rank $i+1$ bundles $\hat V_i$,
  $i\le k$, such that
\[
\nabla \Gamma(\hat V_i) \subset \Omega^1(\hat V_{i+1}), \quad 0\le
i\le k-1\,.
\]
Define
\[
\hat\delta^i = \hat\delta_{i-1} \circ \hat\delta_{i-2}\circ \ldots
\circ \hat\delta_0, \quad 1\le i \le k\,,
\]
and $\hat\delta^0 = \id|_{\hat L}$, where $\hat\delta_i = \pi_{\hat
  V_i} \nabla|_{\hat V_i}$ are the derivatives of the $V_i$.

If $-S$ is adapted to the flag, i.e.,
\[
*\hat\delta_i = -S\hat\delta_i = - \hat\delta_i S
\]
for $ 0 \le i\le k-1$, then the following statements hold:
\begin{enumerate}
\item For all $i=0,\ldots, k$,
 \begin{equation*}
\label{eq:AdeltaQ}
 A\hat\delta^iQ \in H^0(K^{i+2}\Hom_+(V/V_{n-1}, L))
\end{equation*}
is a holomorphic section.  In particular, if $f: S^2\to \HP^n$ is a
Willmore sphere, then $A\hat\delta^iQ= 0$
for all $i=0,\ldots, k$.\\

\item If $A\hat\delta^iQ= 0$ for all $0\le i \le k$, then the image of
  the derivative $\hat\delta_k$ of $V_k$ defines a $\rank k+2$
  subbundle $\hat V_{k+1}\subset V$ provided $\hat\delta_k\not= 0$.
  In this case, the derivative of $V_k$ satisfies
\begin{equation*}
\label{eq:beh} *\hat\delta_{k} = - S\hat\delta_{k} = - \hat\delta_{k} S\,.
\end{equation*}
\end{enumerate}
\end{lemma}
\begin{proof} 
  \quad

\begin{enumerate}
\item Note that $A\hat\delta^iQ\in\Gamma(K^{i+2}\Hom_+(V/V_{n-1},
  L))$ since $S$ is the canonical complex structure of $f$ and thus
  $\ker Q = V_{n-1}$ and $\Im A = L$.
  
  By Lemma \ref{lem:holbundel} the derivatives $\hat\delta_i$ are
  $\partial$--holomorphic.  Since $f$ is a Willmore curve, $A$ and $Q$
  are holomorphic and antiholomorphic respectively, so that \[
  \dbar(A\hat\delta^iQ) = (\dbar A)\hat\delta^iQ +
  A(\partial(\hat\delta^iQ)) = 0\] which shows that $
  A\hat\delta^iQ\in H^0(K^{i+2}\Hom_+(V/V_{n-1}, L))$ is a holomorphic
  section.  Since the degree of a complex holomorphic bundle is given
  by the order of any non-vanishing holomorphic section, we see
\[
\hspace{1cm} \ord(A\hat\delta^iQ)= \deg K^{i+2} + \deg L - \deg V/V_{n-1} \stackrel{\eqref{eq:dualdegree}}=
 (i+2+n)\deg K - \sum_{j=0}^{n-1}\ord\delta_j\,,
\]
provided $A\hat\delta^iQ\not= 0$. If $f: S^2\to\HP^n$ is a Willmore
sphere, then $\deg K <0$ whereas the order of a holomorphic section is
nonnegative. Thus, the above equation cannot hold, and
$A\hat\delta^iQ$ has to vanish.

\item The assumption $A\hat\delta^iQ=0$ for all $i\le k$ implies that
  $\hat V_k\subset\ker A$, i.e. $\hat V_k$ is $A$--stable. But $V_k$
  is $Q$--stable, too, since $\Im Q = \hat L \subset \hat V_k$. In
  view of Lemma \ref{lem:holbundel} it remains to show that $V_k$ is
  also $\partial$--stable.  Lemma \ref{lem:holbundel} then shows that
  $-S$ is adapted and $\hat\delta_k$ is $\partial$--holo\-mor\-phic,
  so that the image of $\hat\delta_k$ defines the smooth bundle $\hat
  V_{k+1}$.
  
  If $k=0$ then $\hat V_k= \hat L = \Im Q$, which is $\partial$ stable
  since $Q$ is antiholomorphic. If $k>0$ then the flag $\hat L \subset
  \ldots \subset \hat V_k$ has $-S$ as an adapted complex structure.
  Hence $\hat\delta_i$ is $\partial$--holomorphic and $\Im
  \hat\delta_{i} = \hat V_i/\hat V_{i-1}$ is $\partial$ stable for
  $0\le i\le k-1$. Thus $\hat V_k$ is again $\partial$ stable.
\end{enumerate}
\end{proof}

Since the backward B\"acklund transform of a Willmore sphere $f:
S^2\to\HP^n$ is a holomorphic curve with respect to the complex
structure $-S$, we can apply the previous Lemma successively as long as
$\hat\delta_k\not=0$ to construct the Frenet flag of $\hat f$. Since
$-S$ is adapted to the flag and $A\delta^kQ=0$, we see that $-S$ is
the canonical complex structure of $\hat f$. Dualizing this result, we
conclude:
\begin{cor}
A B\"acklund transform of a Willmore sphere $f: S^2\to\HP^n$ is a
Frenet curve. Its canonical complex structure is the negative of the
canonical complex structure of $f$.
\end{cor}
\begin{rem}
  The sequence of Willmore spheres obtained by applying successively
  forward (or backward) B\"acklund transformations breaks down after
  at most $n$ steps since the $i^{\text{th}}$ forward B\"acklund
    transform $f_i: S^2\to\HP^{n_i}$ of $f$ maps to $\HP^{n_i}$ with
    $n_i \le n-i$.
\end{rem}

As we have seen before, a Willmore curve $f: M \to\HP^n$ whose
backward B\"acklund transform is a full curve $\hat f: M \to\HP^n$ in
$\HP^n$ and has negative canonical complex structure is a twistor
projection of a holomorphic curve in $\CP^n$.  Conversely, we show
that the B\"acklund transforms of a twistor projection of a
holomorphic curve in $\CP^n$ are Frenet curves:

\begin{cor}
  Let $h: M \to\CP^{2n+1}$ be a holomorphic curve in complex
  projective space such that the twistor projection $f: M \to\HP^n$ of
  $h$ is a Frenet curve.
\begin{enumerate}
\item If $f$ has $Q\not=0$ then the backward B\"acklund transform
  $\hat f: M \to\HP^k$, $k\le n$, of $f$ has $\hat Q = 0$.
\item If the dual curve $f^*$ of $f$ has $A^\dagger\not= 0$ then the
  forward B\"acklund transform $\widetilde{f^*}: M \to\HP^k$, $k\le
  n$, of $f^*$ has $\widetilde{A^\dagger} = 0$.
\end{enumerate}
\end{cor}
\begin{proof}
  Since in the first case $A\hat\delta^kQ = 0$ for all $0\le k \le n$,
  we can construct with Lemma \ref{l:flag} successively flag spaces
  $\hat V_k$ as long as $\hat\delta_k \not=0$.  Let $\hat V$ be the
  $\nabla$ stable bundle so that $\hat f$ is a full curve in $\hat V$.
  By construction, the complex structure $\hat S = -S|_{\hat V}$ is
  adapted, and has $\hat Q = A|_{\hat V} = 0$. In particular, $\hat f$
  is a Frenet curve with canonical complex structure $\hat S$.
  
  The second part is the dual statement of \thetag 1.
\end{proof}

We conclude the paper with a classification result for Willmore
spheres in $\HP^n$.  In the case of a Willmore sphere $f: S^2\to S^4$
the forward and backward B\"acklund transform coincide and give a
point $\infty$ in $\HP^1$ since $A Q=0$. The canonical complex
structure $S$ of $f$ stabilizes $\hat L = \tilde L$. In this case the
canonical complex structure gives the mean curvature congruence of the
conformal map $f: M\to S^4$, and the fact that $S$ stabilizes $\hat L
= \tilde L$ translates to the property that the point $\infty$ lies on
all mean curvature spheres of $f$. Using $\infty$ for a stereographic
projection of $S^4$ to $\R^4$, we see that $f$ becomes a minimal
surface in $\R^4$, see \cite[Sec. \ 11.2]{coimbra} for the details of
this argument. This yields an alternative proof of the result of Ejiri
\cite{ejiri}, see also \cite{montiel}, that a Willmore sphere $f:
S^2\to S^4$ is either a minimal surface in $\R^4$ with planar ends, or
$f$ or its dual curve $f^*$ is a twistor projection of a holomorphic
curve $h: S^2\to \CP^3$.

In the case of Willmore spheres in $\HP^n$ the B\"acklund transforms
are not necessary a constant point but can be Frenet curves in a lower
dimensional $\HP^k$. Since the canonical complex structure of a
B\"acklund transform of a Willmore sphere is $-S$, Proposition
\ref{prop:sphere projection} allows to characterize all Willmore
spheres.

\begin{theorem} 
  Every Willmore sphere $f: S^2\to\HP^n$ has Willmore energy
\[  W(f) \in 4\pi \N,\]
and is given by complex holomorphic data.

More precisely, $f$ is either a minimal surface in $\R^4$ with planar
ends, or $f$ or its dual curve $f^*$ is, at most after projection to a
suitable $\HP^m$, a twistor projection of a holomorphic curve in
complex projective space.

\end{theorem} 

\begin{proof}
If $f: S^2\to\HP^n$ or its dual curve $f^*$ is a twistor projection of
a holomorphic curve in $\CP^{2n+1}$ then $f$ has integer Willmore
energy $W(f)\in 4\pi\N$.

Let $n>1$, $A, Q\not=0$, and assume that one of the B\"acklund
transforms, without loss of generality say $\hat f$, is a full curve
in a rank $n-1$ vector bundle $\hat V$.  Then $\check V = V/\hat V$
has rank 1, so that the projection $\pi: V \to\check V$ only gives a
constant point. However, in this case $\ker A=\hat V$ is a $\nabla$
parallel subbundle of $V$, and $(\ker A)^\perp$ is a constant point.
In other words, $\tilde V$ has rank 1, and the dual curve $f^*$ of $f$
projects to a Willmore curve in $\HP^{n-1}$.

In particular, Proposition \ref{prop:sphere projection} shows for
$n>1$ that at least one of $f$ and $f^*$ projects to the dual curve of
a twistor projection of a holomorphic curve in complex projective
space.  Moreover, $f: S^2\to \HP^n$ has integer Willmore energy
$W(f)\in 4\pi\N$.

The remaining case $n=1$, $A, Q\not=0$, results in minimal spheres in
$\R^4$ with planar ends as discussed above.
\end{proof}




\emph{Acknowledgements.} I would like to thank the members of the
Department of Mathematics and Statistics at UMass for their
hospitality during my stay in Amherst. I'm grateful to the members of
GANG at UMass and SFB 288 at TU Berlin, in particular to Franz Pedit
and to Ulrich Pinkall, for many fruitful and clarifying discussions on
quaternionic holomorphic geometry. I also would like to thank the
referee for many helpful suggestions.

\bibliographystyle{alpha}

\bibliography{doc,doc_dpw}

\newcommand{\etalchar}[1]{$^{#1}$}
\begin{thebibliography}{BFL{\etalchar{+}}02}

\bibitem[BFL{\etalchar{+}}02]{coimbra}
F.~Burstall, D.~Ferus, K.~Leschke, F.~Pedit, and U.~Pinkall.
\newblock {\em Conformal Geometry of Surfaces in ${S}^4$ and Quaternions}.
\newblock Lecture Notes in Mathematics, Springer, Berlin, Heidelberg, 2002.

\bibitem[Bry84]{bryant}
R.~L. Bryant.
\newblock A duality theorem for {W}illmore surfaces.
\newblock {\em J. Diff.\ Geom., Vol. 20}, pages 23--53, 1984.

\bibitem[Eji88]{ejiri}
N.~Ejiri.
\newblock {W}illmore surfaces with a duality in ${S}^n$(1).
\newblock {\em Proc. Lond. Math. Soc., III Ser. 57, No.2}, pages 383--416,
  1988.

\bibitem[FLPP01]{Klassiker}
D.~Ferus, K.~Leschke, F.~Pedit, and U.~Pinkall.
\newblock Quaternionic holomorphic geometry: {P}l\"ucker formula, {D}irac
  eigenvalue estimates and energy estimates of harmonic 2-tori.
\newblock {\em Invent. math., Vol. 146}, pages 507--593, 2001.

\bibitem[FP90]{s4paper}
D.~Ferus and F.~Pedit.
\newblock ${S}^1$-equivariant minimal tori in ${S}^4$ and ${S}^1$-equivariant
  {W}illmore tori in ${S}^3$.
\newblock {\em Math. Z., Vol. 204}, pages 269--282, 1990.

\bibitem[H{\'e}l04]{helein}
F.~H{\'e}lein.
\newblock Removability of singularities of harmonic maps into
  pseudo--{R}iemannian manifolds.
\newblock {\em Ann. Fac. Sci. Toulouse Math., Vol. 1}, pages 45--71, 2004.

\bibitem[LP03]{osculates}
K.~Leschke and F.~Pedit.
\newblock Envelopes and {O}sculates of {W}illmore surfaces.
\newblock arXiv: math.DG/0306150, 2003.

\bibitem[Mon00]{montiel}
S.~Montiel.
\newblock Willmore two-spheres in the four sphere.
\newblock {\em Trans. Amer. Math. Soc., Vol. 352}, pages 4449--4486, 2000.

\bibitem[Pet04]{paule}
P.~Peters.
\newblock {\em Soliton Spheres}.
\newblock PhD thesis, Technische Universit\"at Berlin, 2004.

\bibitem[PP98]{icm}
F.~Pedit and U.~Pinkall.
\newblock Quaternionic analysis on {R}iemann surfaces and differential
  geometry.
\newblock {\em Doc. Math. J. DMV, Extra Volume ICM, Vol. II}, pages 389--400,
  1998.

\bibitem[Rig87]{rigoli}
M.~Rigoli.
\newblock The conformal {G}auss map of submanifolds of the {M}oebius space.
\newblock {\em Ann. Global Anal. Geom., Vol. 5, No. 2}, pages 97--116, 1987.

\bibitem[Wil93]{willmore_book}
T.~J. Willmore.
\newblock {\em Riemannian Geometry}.
\newblock Clarendon Press, Oxford, 1993.

\bibitem[Wol88]{wolfson}
J.G. Wolfson.
\newblock {H}armonic sequences and harmonic maps of surfaces into complex
  {G}rassman manifolds.
\newblock {\em J. Diff.\ Geom., Vol 27}, pages 161--178, 1988.

\end{thebibliography}


\vspace{1cm}

\small
\noindent

\end{document}